\documentclass[10pt]{amsart}

\usepackage{graphicx}

\usepackage{caption}
\usepackage{subcaption}

\newtheorem{thm}{Theorem}[section]

\newtheorem{coro}[thm]{Corollary}
\newtheorem{prop}[thm]{Proposition}

\theoremstyle{definition}

\newtheorem{cond}[thm]{Condition}

\def\ss{\mathbb{S}^n}

\def\refl{\text{Refl}}
\def\ant{\text{Ant}}

\def\ZZ{\mathbb{Z}}
\def\CC{\mathbb{C}}
\def\RR{\mathbb{R}}

\title{Isospectral configurations in Euclidean and Hyperbolic Geometry}
\author[Masai]{Hidetoshi Masai}
\address{Department of Mathematics Tokyo Institute of Technology 2-12-1, Ookayama, Meguro-ku, Tokyo. 152-8551. Japan}
\email{masai at math.titech.ac.jp}
\author[McShane]{Greg McShane}
\address{Institut Fourier 100 rue des maths, BP 74, 38402 St Martin d'H\`eres cedex, France}
\email{mcshane@univ-grenoble-alpes.fr}
\thanks{This work was supported by PERSYVAL-LAB/Equipe Action TOFU and JSPS KAKENHI Grant Number 19K14525 and 21H04428.}
\begin{document}
\maketitle

\section*{Abstract} 
A number of questions related to the length spectrum of surfaces are discussed and in particular the existence of pairs of surfaces which though not isometric are isospectral.
Here by isospectral we mean that a pair of bodies have the same distribution of chord lengths.
In the Euclidean setting, we study isospectral convex dodecagons found by Mallows and Clark in the  1970’s.
Starting from their idea, we give constructions for isospectral pairs of hyperbolic surfaces that have no common cover.
Since the work of Mallows and Clark is probably unfamiliar to readers with a background in topology/hyperbolic geometry
we include expository material on other related topics about the distribution of chord lengths.

\section{Introduction}
\subsection{Chord lengths}
The distribution of lengths of chords of convex bodies 
is a fundamental problem in integral geometry
and has applications to  tomography and x-ray crystallography.
The basic question is: 
does the distribution of chord lengths determine
the convex body up to isometry.
By cleverly partitioning 
the sides of a regular octagon
into two  subsets $X$ and $Y$ 
Mallows and Clark \cite{MC} constructed
a pair of convex dodecagons
which are not congruent but have 
the same distribution of chord lengths.
They obtain  these dodecahedra 
by suitably capping off
each of the sides of $X$ (resp. $Y$)
with a triangle.

What is important in the Mallows Clark construction,
as was observed by R. Garcia-Pelayo \cite{GP},
is that the configuration of sides 
in $X$ and $Y$ are associated 
to \textit{non congruent homometric pairs} (NCHP) in the cyclic
graph  $C_8$.
He further observes that, 
after applying an isometry of $\mathbb{R}^2$
the sets $X$ and $Y$ can be chosen to be complementary.
Garcia-Pelayo goes on to show that
any pair of complementary sets,
of a  finite vertex transitive graph 
are in fact homometric
thus generalizing a theorem of A. L. Patterson from the 1940s.
Here by \textit{pair of complementary sets}
we mean a pair of subsets of vertices of the same cardinal
whose reunion is the set of all vertices.

\begin{figure}[hb]
\centering
  \includegraphics[scale=.4]{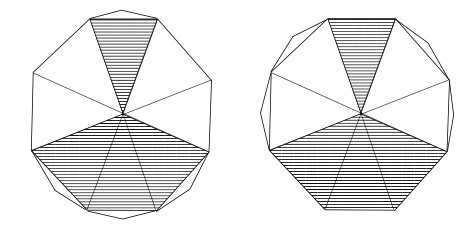} 
  \caption{Mallows Clark pair.}
  \label{mc pair}
\end{figure}

\subsection{Spanning rulers}
In fact, NCHP have been extensively studied in relation
to spanning rulers.
 A \textit{ruler} is any finite subset of the positive integers
 containing $0$ and the points 
 of the set are called  \textit{marks}.
 A \textit{spanning ruler} is a ruler
such that the distances between 
(distinct) points are distinct.
More informally  a \textit{regular ruler}
consists of the points$\{0,1,2\ldots n \}$
so that the distance $1$ can be measured in 
$n-1$ ways (e.g $2-1$, $3-2$ etc.)
whilst on a spanning ruler  
each distance that can be measured
can be done so  in exactly one way.
An example of a spanning ruler is 
$0,1,4, 6$ 
which  measures
 all the distances $1,2,3,4,5,6$ just once.
There are diverse applications of  spanning rulers
ranging from  in  radio astronomy to x-ray diffraction crystallography.

In 1939 Sophie Piccard claimed that 
for any spanning ruler the set of distances it measures
determines it up to isometry
or that, equivalently,  
any ruler from an NCHP  is never a spanning ruler.
In the 70's G.  Bloom \cite{Bloom} 
found an NCHP  which  is a counterexample to her claim, namely
the pair
$ \{0,1,4,10,12,17\}$ and $\{0,1,8,11,13,17\}$.
Subsequently, 
Bloom-Golomb \cite{Bloom-Golomb}
 found an  infinite family of counterexamples 
 each with six marks containing this pair.

\subsection{Structure and results}

Our principal motivation was an attempt to find pairs of hyperbolic
surfaces that had the same chord distribution which are not
isometric and perhaps not even commensurable (i.e. they have no
common cover of finite degree.
Thus, with the exception of Sections \ref{section:hyperbolic} and
\ref{section:projecting} which have some geometric/topological
content most of this text is expository.

We begin (Section \ref{section:auto-correlation} and
\ref{section:swapping}) by reformulating Garcia-Pelayo's result in
terms of auto-correlation functions which will allow us to present a
short proof of them. In Section \ref{section:pairs of subsets} we discuss a
particular family of NCHP using this approach. Further we
construct  NCHP  which are counterexamples to the
questions asked by  Garcia-Pelayo.
In particular we find NCHP pairs for $n$-polygons for $n>15$ odd
using the topological ideas introduced in Section
\ref{section:projecting}.

In Section \ref{section:hyperbolic} we adapt the construction of Mallows Clark in the hyperbolic setting
yielding examples of non congruent  ideal polygons with the same
distribution of lengths of chords.
We also construct other pairs of
hyperbolic surfaces with the same distribution of lengths of chords.
All our examples are of rather simple topological type being
homeomorphic to annuli. 
It seems difficult to construct non isometric examples of pairs of more general surfaces
with the same distribution of chord lengths; see for example the
discussion in \cite{Masai-McShane}.
Nonetheless even such simple surfaces have appeared in several
papers in the context of hyperbolic geometry and identities for the
Roger's Dilogarithm \cite{McShane-identities,Bridgeman, Tan_et_al}.

We end the paper with two appendices with which we hope to tempt the
reader to engage further in the problem of construction and
classification of NCHP. The first describes the Mallows-Clark
construction for NCHP pairs in geometric terms and the second
Golomb's algebraic approach using what amounts essentially to
factorisation of the auto-correlation function.

\section{Generating and auto-correlation
functions}\label{section:auto-correlation}
Let $G$ be a connected  graph with vertices $V(G)$.
 There is a natural \textit{path metric} on $G$
 with distance function $d: V(G) \times V(G) \rightarrow \RR_+$ 
such that  for  adjacent vertices 
$v_1,v_2 \in V(G)$  the distance satisfies $d(v_1,v_2) = 1$.
It is easy to see that any automorphism of $G$ 
is an isometry of this metric.
To each  vertex $v\in V(G)$ we associate
a \textit{generating series} namely
$$\sigma_{v}(t)  = \sum_{u \in V(G)} t^{d(v,u)}.$$
If $(X_i)_i$ is a decomposition of $V(G)$ that is 
$$ V(G) = \bigsqcup_i X_i,$$
then $\sigma_v$ decomposes as 
$$\sigma_{v}(t)  = \sum_i \left( \sum_{u \in X_i} t^{d(v,u)} \right).$$
More generally define the \textit{auto-correlation 
function} of $X \subset V(G)$
to be 
$$\sigma_X(t) = 
\frac12 \sum_{v \in X} \sigma_v(t) = 
\frac12 \sum_{(v,u) \in X^2} t^{d(v,u)}  .$$
Let $G$ be a connected  graph 
and $T$ a fixed point free automorphism of $G$
then there is an obvious relation between these series:
If $X\subset V(G)$ is such that
$V(G) = X \sqcup T(X) $
then by letting $\overline{X}:=V(G)\setminus X=T(X)$, we have
$$\sigma_X = \sigma_{\overline{X}}.$$
\section{Swapping in vertex transitive graphs}
\label{section:swapping}
Recall that a graph $G$  is \textit{vertex transitive}
iff its automorphism group acts transitively on the vertices $V(G)$.
Note that a vertex transitive graph must be \textit{regular},
that is every vertex has the same valence.
Let $G$ be a vertex transitive graph.
Since any automorphism is an isometry of the natural path metric
one has, for all $u,v \in V(G)$,
 $$ \sigma_u (t)  =  \sigma_v (t). $$
 Let $X$ be a subset of $V(G)$ and 
  ${\overline{X}}$ denote the complement of $X$.
A  set $Y$ is obtained from $X$   by \textit{swapping} $u \in X$ for  $ v \in \overline{X}$ iff
$$
X = Z  \cup \{u\},\, Y = Z  \cup  \{v\},
$$
where $Z = X \cap Y$.
Note that $ \overline{Y}$  is obtained from $\overline{X}$
by swapping $v$ for $u$ since
$$\overline{X} = Z' \cup \{v\},\, \overline{Y} = Z' \cup \{u\},$$
where   $Z' = \overline{X}  \cap \overline{Y}$.
If $X_0$ and $X_1$ are a pair of finite subsets with the same 
number of elements then one can transform 
$X_0$ into $X_1$ by a finite number of swaps.
The minimal number of swaps is called the \textit{Hamming distance}.

\begin{thm}\label{difference invariant by swapping}
Let $G$ be a connected finite vertex transitive  graph 
and $X \subset V(G)$ and  
 ${\overline{X}}$ denote the complement of $X$.
Then for any subset $Y$ obtained from $X$ by a finite number of swaps 
$$\sigma_X - \sigma_{\overline{X}}  = \sigma_Y - \sigma_{\overline{Y}}.$$

\end{thm}

\begin{proof}
The proof is by induction on the number of \textit{swaps.}

It is easy to see that it suffices to prove the theorem for
a single swap. We present a proof 
which is  completely formal as follows.
There is an obvious decomposition
  of the spectral functions:
\begin{eqnarray*}
 \sigma_X = \sigma_{Z}  
 + \sum_{x \in Z}  t^{d(u,x)} + \frac12 \\
  \sigma_{\overline{X} } =
   \sigma_{Z'}  
  + \sum_{x \in \overline{Z'}}  t^{d(v,x)}  + \frac12
\end{eqnarray*}  
likewise
\begin{eqnarray*}
 \sigma_Y = \sigma_{Z}  
  +\sum_{x \in Z}  t^{d(v,x)}  + \frac12\\
  \sigma_{\overline{Y} } =
   \sigma_{Z'}  
 + \sum_{x \in \overline{Z'}}   t^{d(u,x)} + \frac12 
\end{eqnarray*}  

Thus

\begin{eqnarray}
\sigma_{X} - \sigma_{\overline{X} } 
=  
\sigma_{Z } 
- \sigma_{Z'} + 
\left( \sum_{x \in Z}  t^{d(u,x)} 
- 
 \sum_{x \in Z'}  t^{d(v,x)} \right)  \label{aa}\\ 
\sigma_{Y} - \sigma_{\overline{Y} } 
=  
\sigma_{Z } 
- \sigma_{Z'} + 
\left( \sum_{x \in Z}  t^{d(v,x)} 
- 
 \sum_{x \in Z'}  t^{d(u,x)} \right).\label{bb}
\end{eqnarray} 
The generating functions also decompose as
\begin{eqnarray*}
 \sigma_u (t)  =  \sum_{x \in Z}  t^{d(u,x)}
 + \sum_{x \in Z'}  t^{d(u,x)}  +  t^{d(u,v)} \\
  \sigma_v (t)  =   \sum_{x \in Z}  t^{d(v,x)}
 + \sum_{x \in  Z'}  t^{d(v,x)}  +  t^{d(u,v)} \\
\end{eqnarray*}  
By hypothesis $ \sigma_u (t)  =  \sigma_v (t) $ 
so subtracting these expressions 
one sees that the expressions in parentheses 
in equations (\ref{aa}) and (\ref{bb}) are equal.
\end{proof}

\begin{coro}Let $G$ be a connected finite vertex transitive graph.
If $\sigma_X =  \sigma_Y$ then 
$\sigma_{\overline{X} } = \sigma_{\overline{Y} }$.
\end{coro}

\begin{proof}
By the theorem 
\begin{eqnarray*}
\sigma_X - \sigma_{\overline{X}}  &=& \sigma_Y - \sigma_{\overline{Y}}  \\
&=&  \sigma_X - \sigma_{\overline{Y}}
\end{eqnarray*}
so the result follows.
\end{proof}

\begin{coro}\label{comp pairs}
Let $G$ be a  connected finite vertex transitive graph with $2n$ vertices.
Suppose $G$ has a fixed point free automorphism $T$ of even order.
Then for any subset of  $X$ of $n$ vertices
$$\sigma_X =  \sigma_{\overline{X} }.$$
\end{coro}

\begin{proof}
Since $T$ is of even order
there  exists a subset $X_0 \subset V(G)$,
such that  $T(X_0) = \overline{X_0}$.
Further $T$ being a bijection,
the sets $T(X_0)$ and $X_0$ have the same 
number of elements,
so in particular $X_0$ must have  exactly $n$ elements.
Since $T$ is an isometry of the path metric
$\sigma_{X_0} =  \sigma_{\overline{X}_0 }.$

Since $V(G)$ is finite  
one can transform any other subset $X$ of $n$ element
into $X_0$
by a finite number of swaps so it follows  from Theorem 
 \ref{difference invariant by swapping}
that 
$$0 = \sigma_{X_0} -   \sigma_{\overline{X}_0 }
= \sigma_{X} -   \sigma_{\overline{X} }.$$
\end{proof}

\section{Pairs of subset}\label{section:pairs of subsets}

Using Corollary \ref{comp pairs} it is easy 
to find pairs of  sets $X$ and $Y$ such that 
$\sigma_X = \sigma_Y$ but which are not congruent.
Here we say that  $X$ and  $Y$
are \textit{congruent} 
if there is  an automorphism of $G$ mapping $X$ onto $Y$.
In fact, it is also easy to  find 
pairs of  sets $X$ and $Y$ such that 
$\sigma_X = \sigma_Y$ 
and  $Y$ is congruent to 
neither  $X$ nor $\overline{X}$.

Let $G$ be a cycle of $2n$ vertices and   $X$ 
a subset of $n$ vertices such that $\overline{X}$
is not congruent to $X$.
Let $2G$ denote the graph obtained from $G$
by subdividing each edge, 
it is a  cycle of $4n$ vertices.
There is a natural inclusion $\iota:V(G) \hookrightarrow V(2G)$
and, since $d(\iota(u),\iota(v)) = 2 d(u,v)$, one has
$$\sigma_{\iota(X)}(t) = \sigma_{X}(t^2) ,\,\
\sigma_{\iota(\overline{X} )}(t) = \sigma_{\overline{X}}(t^2).$$
But $\iota(\overline{X}) $ has $n$ vertices and so 
is not congruent to the complement of $\iota(X)$ 
which has $3n$.
\subsection{Homometric configurations on the line}
For an $n$-point configuration $X = \{x_{0},\dots,x_{n-1}\}$
the distance spectrum $\mathcal{D}(X)$ is defined as
the set  $$\{|x_{i}-x_{j}|\mid 0\leq i<j\leq n-1\},$$
counted with multiplicities.
It is convenient to label the points so that
$x_{i}<x_{i+1}$ for all $0\leq i\leq n-1$.
Two configurations are equivalent if they are related by a finite sequence of translations, dilations, and reflections of the line.
The following proposition provides us with 
a recursively defined
 infinite family of inequivalent homometric configurations.
\begin{prop}
For every $k\geq 12$, 
the sets $X_{k}:=\{0,1,5,6,7,9,10...,k-4,k-3,k\}$ and
$Y_{k}:=\{0,4,5,6,9,10,...,k-4,k-3,k-1,k\}$
are $(k-5)$-point inequivalent homometric configurations.
\end{prop}
\proof
Note that 
\begin{align*}
X_{k+1} &= X_{k}\setminus \{k\} \cup \{k-2,k+1\}\\
Y_{k+1} &= Y_{k}\setminus\{k-1\} \cup \{k-2,k+1\}.
\end{align*}
Hence to compute $\mathcal{D}(X_{k+1})\setminus \mathcal{D}(X_{k})$ and 
$\mathcal{D}(Y_{k+1})\setminus \mathcal{D}(Y_{k})$, it suffices to deal with three points in each case.
Thus it can be seen that both 
$\mathcal{D}(X_{k+1})\setminus \mathcal{D}(X_{k})$
 and $\mathcal{D}(Y_{k+1})\setminus \mathcal{D}(Y_{k})$
consist of the points
$$\{1,2,\dots,k-11,k-9,k-8,k-8,k-4,k-3,k-2,k+1\} - \{k-1\}.$$ \qed
\subsection{Configurations of five  or less points}

\begin{figure}
\centering
  \includegraphics[scale =0.3]{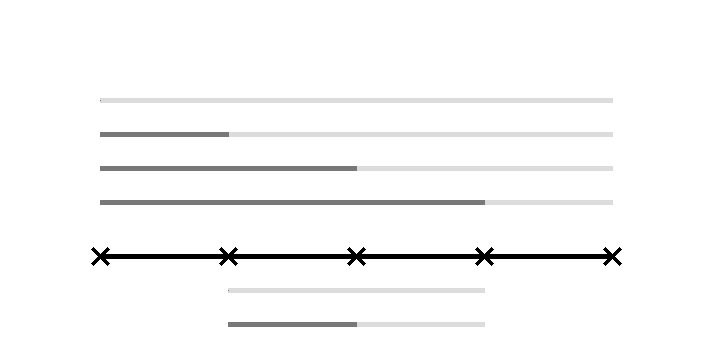}
  \caption{5 points configuration.}
  \label{fig.5pts}
\end{figure}

It is  known that 
there are no  homometric pairs with  less than six points.
We give a short proof of this fact for  completeness.
Let $X = \{v_{1}, \dots, v_{5}\}$ 
be five points on the line arranged in order according to the index
and note that 
there are exactly  10 elements in the
 set of distances  $\mathcal{D}(X)$.
Note that $|v_5 - v_1|$ 
is the maximum of  $\mathcal{D}(X)$
and that,
with  the exception of
$|v_5 - v_1|$ and $|v_4 - v_2|$, 
the distances can be paired off
(see  Figure \ref{fig.5pts}) 
so that the  sum of a pair is either 
$|v_5 - v_1|$ or $|v_4 - v_2|$
thus:
\begin{equation}\label{sum of lengths}
2|v_{4}-v_{2}| +  4|v_{5}-v_{1}| = \sum_{l\in\mathcal{D}(X)}l 
\end{equation}
The second largest distance in $\mathcal{D}(X)$
is either $v_{5} - v_{2}$ or $v_{4} - v_{1}$.
Since we know the  value of  $v_{4}-v_{2}$ 
from  formula  (\ref{sum of lengths}) above 
we can determine $v_2$ and $v_4$
and the the 6 distances $|v_i - v_j|,\,  i >  j , i \neq 4, j \neq 3$.
After removing these distances
there are 2 pairs of distances
involving $v_3$ and we can determine 
$v_3$ from the equations:
\begin{eqnarray*}
|v_3 - v_2| + |v_3 - v_4 | & = &  |v_4 - v_2| \\
|v_3 - v_1| + |v_3 - v_5  | & = &  |v_5 - v_1|.
\end{eqnarray*}

\section{Ortho spectrum} \label{section:hyperbolic}
We consider the hyperbolic plane $\mathbb{H}^{2}$.
Let $\mathcal{G}$ be a collection of mutually disjoint (possibly asymptotic) geodesics.
Two such collections are {\em congruent} if there is an isometry which maps one to the other set-wise.
An {\em ortho geodesic} is a geodesic each of two endpoints is orthogonal to an element of $\mathcal{G}$.
The {\em ortho spectrum} $\mathcal{O}(\mathcal{G})$ is the set of the lengths of ortho geodesics counted with multiplicity.
For a hyperbolic surface with totally geodesic boundary, we get a collection of geodesics as the boundary of the universal covering, and ortho spectrum is defined to be the one for those geodesics.
In this section, by using Corollary \ref{comp pairs}, we construct examples of incongruent collections of geodesics with the same ortho spectrum.
Some of those examples also give examples of non-isometric hyperbolic surfaces with the same ortho spectrum.
In \cite{Masai-McShane}, we further discuss systoles and ortho spectrum rigidity.
First, let $g$ be a geodesic, and $\gamma$ be an isometry which satisfies the following condition.

\begin{cond}\label{cond.isom}
For all $i\in\mathbb{Z}$,
\begin{itemize}
\item we have either $g\cap \gamma^{i} g = \emptyset$ or $g = \gamma^{i}g$, and
\item $\gamma^{i}g$'s are in the same component of the complement of $g$.
\end{itemize}
\end{cond}

\begin{figure}
\centering
\begin{minipage}{.25\textwidth}
  \centering
  \includegraphics[width=1\linewidth]{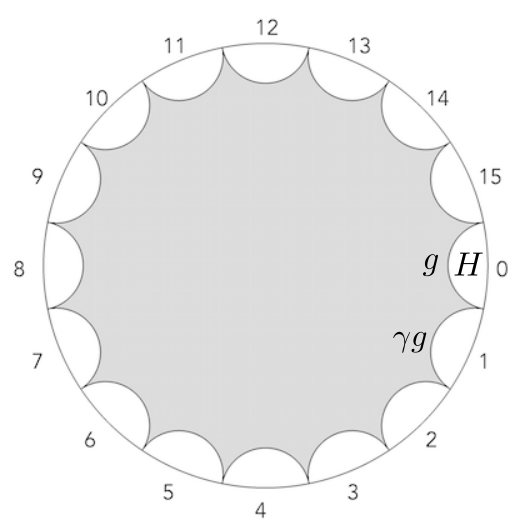}
\end{minipage}%
\begin{minipage}{.35\textwidth}
  \centering
  \includegraphics[width=1\linewidth]{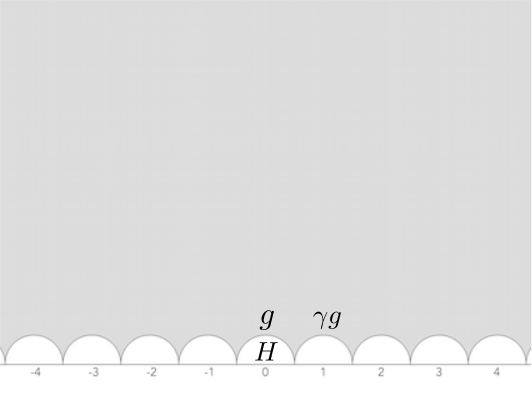}
\end{minipage}
\begin{minipage}{.35\textwidth}
  \centering
  \includegraphics[width=1\linewidth]{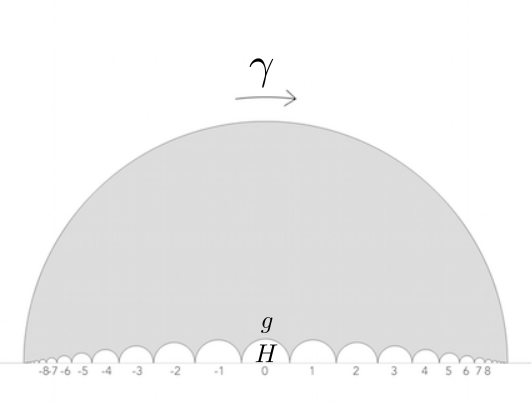}
\end{minipage}
\caption{(left) elliptic, (center) parabolic, (right) hyperbolic.}
\label{fig.translations}
\end{figure}

Let $H$ be the half space with respect to $g$ that does not contain $\gamma^{i}g$'s.
Let $\mathcal{G}_{0}$ be a collection of disjoint geodesics in $H$.
Then put $\mathcal{G}_{i}:=\gamma^{i}(\mathcal{G}_{0})$ for $i\in \mathbb{Z}$.

For two disjoint collections of geodesics $\mathcal{H}_{1}$ and $\mathcal{H}_{2}$,
let $\mathcal{O}(\mathcal{H}_{1},\mathcal{H}_{2})$ be the lengths of all ortho geodesics connecting $\mathcal{H}_{1}$ and $\mathcal{H}_{2}$, counted with multiplicity.
The following proposition is immediate from the construction.
\begin{prop}\label{pair spec}
$\mathcal{O}(\mathcal{G}_{i_{1}}, \mathcal{G}_{i_{2}}) = \mathcal{O}(\mathcal{G}_{j_{1}}, \mathcal{G}_{j_{2}})$ 
whenever $|i_{1}-i_{2}| = |j_{1}-j_{2}|$.
\end{prop}
Let $C_{2n}$ be the cycle graph of length $2n$ and 
$\widetilde{C_{2n}}$ a universal covering of $C_{2n}$.
We choose a generator $t$ of covering transformation and label $V(\widetilde{C_{2n}})$ so that for each $i\in\mathbb{Z}$,
\begin{itemize}
\item $t(v_{i}) = v_{i+2n}$, and 
\item $v_{i-1}$ and $v_{i}$ is connected by an edge.
\end{itemize}
Let $X'$ be an $n$-point subset of $V(C_{2n})$ and $X$ the pre-image to $\widetilde{C_{2n}}$ of $X'$.

We define two collections of geodesics by
\begin{eqnarray*}
\mathcal{G}_{X}:= \bigcup_{v_{i}\in X}\mathcal{G}_{i},\\
\mathcal{G}_{\overline{X}}:= \bigcup_{v_{i}\in \overline{X}}\mathcal{G}_{i}.\\
\end{eqnarray*}
Then we have
\begin{thm}
Two collections $\mathcal{G}_{X}$ and $\mathcal{G}_{\overline{X}}$ have the same ortho spectrum, that is
$\mathcal{O}(\mathcal{G}_{X}) = \mathcal{O}(\mathcal{G}_{\overline{X}})$.
\end{thm}
\proof 
Let $I_{X}:=\{(i,j)\mid i<j , v_{i}, v_{j}\in X$\}.
If $\gamma^{m}g = g$ for some $m$, then we only consider $0\leq i<j<m$, see also Corollary \ref{coro.elliptic} below.
The ortho spectrum can be decomposed as
$$\mathcal{O}(\mathcal{G}_{X}) = 
\bigcup_{v_{i}\in X} \mathcal{O}(\mathcal{G}_{i}) \cup
\bigcup_{(i,j)\in I_{X}} \mathcal{O}(\mathcal{G}_{i},\mathcal{G}_{j}).$$
Hence by Proposition \ref{pair spec}, 
we see that $\mathcal{O}(\mathcal{G}_{X})$ is determined by the spectral series of $X$.
Note that Proposition \ref{pair spec} includes the case where $|i-j|=0$.
By Corollary \ref{comp pairs}, $X$ and $\overline{X}$ have the same spectral series and hence we have $\mathcal{O}(\mathcal{G}_{X}) = \mathcal{O}(\mathcal{G}_{\overline{X}})$. \qed
\begin{coro}
There are incongruent collections of geodesics with the same ortho spectrum.
\end{coro}
\proof
Any incongruent pair $X'$ and $\overline{X'}$ of $n$-point subsets of $V(C_{2n})$ gives a incongruent pair of geodesics $\mathcal{G}_{X}$ and $\mathcal{G}_{\overline{X}}$ with the same ortho spectra.
Incongruent pairs of $n$-point subsets of $V(C_{2n})$ are given in \cite{MC, GP}.
\qed

According to the type of the isometry $\gamma$, we have three different kinds of examples.
Let $H_{X}$ (resp. $H_{\overline{X}}$) be the connected component of $\mathbb{H}^{2}\setminus \mathcal{G}_{X}$
(resp. $\mathbb{H}^{2}\setminus \mathcal{G}_{\overline{X}}$) whose boundary contains whole $\mathcal{G}_{X}$ 
(resp. $\mathcal{G}_{\overline{X}}$).

\begin{coro}\label{coro.elliptic}
There are non-isometric hyperbolic orbifolds with the same ortho spectrum.
\end{coro}
\proof
We suppose that $\gamma$ is elliptic.
Then Condition \ref{cond.isom} implies that $\gamma$ is of finite order.
Let $m$ be the smallest integer so that $\gamma^{m}g = g$.
We also consider the case where $\gamma$ is identity, or $m=0$.
We assume that $2n$ divides $m$.
Both $H_{X}$ and $H_{\overline{X}}$ are symmetric with respect to $\gamma^{2n}$.
Hence $H_{X}/\langle \gamma^{2n}\rangle$ and 
$H_{\overline{X}}/\langle \gamma^{2n}\rangle$ are hyperbolic orbifolds.
It is easy to see that if $X$ and $\overline{X}$ are incongruent, then $H_{X}/\langle \gamma^{2n}\rangle$ and 
$H_{\overline{X}}/\langle \gamma^{2n}\rangle$ are non-isometric.
\begin{coro}
There are incommensurable hyperbolic surfaces with the same ortho spectrum.
\end{coro}
\proof
We give two kinds of examples.
First, suppose $\gamma$ is parabolic.
We assume that the endpoints at infinity of $g$ are not parabolic fixed point of $\gamma$.
Then for any $m$ with $2n|m$, $H_{X}/\langle \gamma^{m}\rangle$ and 
$H_{\overline{X}}/\langle \gamma^{m}\rangle$ are punctured surfaces with the same ortho spectrum.
They are non-isometric whenever $X$ and $\overline{X}$ are incongruent. 

If $\gamma$ is hyperbolic, 
we assume that the geodesic axis $h$ does not intersect with $g$.
Then let $H'$ be the half space with respect to $h$ which contains $g$.
Then for any $m$ with $2n|m$, $(H_{X}\cap H')/\langle \gamma^{m}\rangle$ and 
$(H_{\overline{X}}\cap H')/\langle \gamma^{m}\rangle$ are hyperbolic surfaces with one totally geodesic boundary.
They have the same ortho spectrum, and again they are non-isometric whenever $X$ and $\overline{X}$ are incongruent.

Since the sets of consecutive distances of the Mallows-Clark pair are different,
we see that they are incommensurable. 
Hence the pairs of hyperbolic surfaces we constructed above are also incommensurable.
\qed

\section{Projecting from $\ZZ$ onto $C_n$}\label{section:projecting}

In this section we show how to get examples 
of non congruent homometric pairs of subsets
of the cyclic graph on $n$ vertices using 
covering maps.
In particular we obtain subsets which are 
not pairs of complementary sets as in Section \ref{section:swapping}.
Further we find NCHP pairs for $n$-polygons for $n>15$ odd
using the topological ideas introduced in the previous section.

\subsection{Distances and covering maps}
Let $C_n$ be the cyclic graph on  $n$ vertices  and
$C_\infty$ the connected  infinite graph all
 of whose vertices are valence 2.
 We can identify $C_\infty$ with the universal cover of $C_n$
 and there is a projection map $\pi :C_\infty \rightarrow C_n.$
One can view these graphs as Cayley graphs of  
(cyclic) groups in an obvious way:
 $C_\infty$ is  the Cayley graph of   $(\ZZ,+)$ 
with respect to  the generator $1$
 and $C_n$ of $(\ZZ/n\ZZ,+)$.
 In this way  $\ZZ$  is naturally  identified
 with  the vertices  of   $C_\infty$   
 and $\ZZ/n\ZZ$ with those of $C_n$.
 In particular we can identify 
 the vertices of the latter with $0,1,\ldots n-1$
 in the usual way.

Let $x_i,y_i \in C_\infty$, 
if $d(x_1,y_1) = d(x_2,y_2) $ then 
$d(\pi(x_1), \pi(y_1) ) = d(\pi(x_2) , \pi(y_2) )$.
Since $C_\infty$ is a vertex transitive 
it suffices to fix $x$ and 
to check that there is a function $F$ such that 
for every vertex $y$ of $C_\infty$ 
$$d( \pi(x) , \pi(y) )  = F(d( x, y )).$$
One checks that $F$ is the absolute value of the 
symmetric remainder for the division by $n$.
Recall that  \textit{symmetric remainder}
is the integer $r$ such that $ p = qn +r$
and $|r| <  n/2$.

\subsection{Images of homometric pairs under covering maps}
An $n$-point  configuration $X$  
is a map $\iota_X : \{ 1,2,\ldots n \}  \rightarrow \ZZ$.
It is convenient to work with configurations 
as they correspond to subsets of points counted
with multiplicities.
We say that a point $x$  of $X$ is a 
\textit{multi point}
iff $\{x\}$ is strictly contained in 
 $\iota_X^{-1} (\{ \iota_X(x) \})$.
 Two configurations $\iota_{X}$ and $\iota_{Y}$ are isometric iff
 there exists $g$ an isometry of $\ZZ$ 
 such that $g\circ\iota_X = \iota_Y$.
 Any subset of $\ZZ$ can be viewed as a configuration
 with no multi points.

If $\tilde{X} ,\tilde{Y}$ is  a pair
of non congruent  homometric subsets  of $\ZZ$
then, by the above discussion,
 under the projection map $\pi$ 
the resulting configurations  $\pi(\tilde{X})$ and $\pi(\tilde{Y})$  are homometric.
The examples we obtain by this method fall into three classes:
\begin{enumerate}
\item If $n$ is greater than twice the diameter of 
$\tilde{X}$ then 
the restriction of the projection $\pi$ to $\tilde{X}$ (resp. 
$\tilde{Y}$)
 is injective and, moreover, 
 this restriction maps $\tilde{X}$ isometrically onto 
 $\pi(\tilde{X})$.
 Thus in this case we obtain a pair of non congruent 
homometric subsets of $C_n$;
\item If $n$ is greater than  the diameter of $\tilde{X}$
 but less than  twice the diameter of $\tilde{X}$
 then
the restriction of the projection $\pi$ to $\tilde{X}$ (resp. $\tilde{Y}$)
 is injective  but   $\tilde{X}$  is not mapped
  isometrically onto $\pi(\tilde{X})$.
  Thus we obtain a pair of (possibly congruent)
  homometric sets of $C_n$.
\item If $n$ is  less than   the diameter of $\tilde{X}$ then  
the restriction of $\pi$ is not a priori  injective
so we obtain (possibly congruent) homometric 
 configurations.
 \end{enumerate}

\subsection{An example}
The sets 
$\tilde{X} = \{ 0, 2, 5, 14, 18, 25 \}$ and 
$\tilde{Y} = \{ 0, 2, 13, 16, 20, 25 \}$ are non congruent  homometric subsets of $\ZZ$ and under the projection
to $\ZZ/n\ZZ$ for $n = 11,12,13,14$ 
they give configurations  each with a single double point.
For $n = 15$ one obtains a pair of sets
$X = \{ 0, 2,3, 5, 10, 14 \},  Y =  \{ 0, 1, 2, 5, 10, 13 \}$.
These sets are not congruent 
in $C_{15}$ since the  distances between consecutive points are different: these are respectively 
$2,1,2,3,4,1$ and  $1,1,3,5,3,2$.

Thus we obtain a table of pairs of NCHP:

\begin{center}
\begin{tabular}{|c|l|l|}
\hline  15 & 0, 2, 3, 5, 10, 14 & 0, 1, 2, 5, 10, 13 \\
\hline  17 & 0, 1, 2, 5, 8, 14 & 0, 2, 3, 8, 13, 16 \\
\hline  19 & 0, 2, 5, 6, 14, 18 & 0, 1, 2, 6, 13, 16 \\
\hline  21 & 0, 2, 4, 5, 14, 18 & 0, 2, 4, 13, 16, 20 \\
\hline  22 & 0, 2, 3, 5, 14, 18 & 0, 2, 3, 13, 16, 20 \\
\hline  24 & 0, 1, 2, 5, 14, 18 & 0, 1, 2, 13, 16, 20 \\
\hline 
\end{tabular} 

\end{center}
For $C_{21}$ the two sets we obtain are 
$0, 2, 4, 5, 14, 18$ and  $0, 2, 4, 13, 16, 20$
for which  the distances between consecutive points are 
respectively 2,2,1,9,4,3 and 2,2,9,3,4,1.
So these sets are not isometric even though
(setwise)  these distances between consecutive points  are the same.

\section{Remarks, further questions}

As indicated in the introduction, 
in Section \ref{section:hyperbolic} 
we succeeded in constructing  pairs of
hyperbolic surfaces with the same distribution of lengths of chords.
Our examples are of rather simple topological type being
homeomorphic to annuli and in particular what Casson and Bleier
\cite{Casson}
refer to as crowns - that is there is a boundary component that
consists of a single closed geodesic and the other boundary is a
union of complete geodesics asymptotic in a finite numberof spikes.

\noindent \textbf{Question:} Is it possible to construct a pair of non
isometric surfaces, each of  which consists of a pair of crowns
identified along the boundary components that are closed geodesics,
but which have a common orthospectrum?

\section{Appendix: Structure of homometric pairs}\label{appendix: structure of pairs}

In this section we study the underlying structure of the
Mallows-Clark pair decomposing it into 
three sets $X$, $Y$ and $P_\pm$. 
Our decomposition will allow us to show that:
\begin{itemize}
\item the Mallows-Clark pair lies in a one parameter 
family of NCHP of the circle.
\item to construct, for each $n\geq 4$,  NCHP 
in the $n$-sphere $\ss$.
\end{itemize}

\subsection{Structure of the Mallows-Clark pair}
\label{structure}

Each of the two  configurations  $C_1,C_2$
of the  Mallows-Clark pair consists of four points and
for each of the configurations
exactly two of these points lie on a diameter of the circle.
We can apply a rotation so that 
$C_1 \cap C_2$ consists of 
exactly three points: 
a pair on the diameter  $x_1,x_2$
and another  point $y$.
The set  $C_1 \setminus C_2$ consists of a single  point $p_1$ 
likewise  $C_2 \setminus C_1$ is just $p_2$.
One verifies that $p_1,p_2$ are also on a diameter.
Thus
\begin{eqnarray*}
C_1 \cup C_2& =& X \sqcup Y \sqcup P \\
C_1 & =& X \sqcup Y \sqcup \{ p_1\}\\
C_2 & =& X \sqcup Y \sqcup \{ p_2\}
\end{eqnarray*}
where $X = \{ x_1,x_2\}$, $Y = \{ y\}$
and $P = \{p_1,p_2\}$.

Let $\ant$ denote the antipodal map of the circle
and $\refl$ the unique orientation reversing map
that has $y$ as its fixed point.
Both of these maps are isometries of the circle.
Then:
\begin{enumerate}
\item $X$ is invariant under the antipodal map $\ant$;
\item $Y$ is invariant under $\refl$;
\item $P$ is invariant under both $\ant$ and $\refl$
and each of these maps exchange $p_1$ and $p_2$
that is
\begin{equation}
\ant(p_1)  = p_2,\, \refl(p_1)  = p_2.
\end{equation}
\end{enumerate}

\begin{center}
\includegraphics[scale=.7]{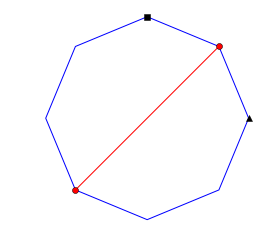} 
	\captionof{figure}{One of the configurations of a
	Mallows-Clark pair. 
	The red dots are antipodeal pairs. 
	The black triangles is a point  invariant
	under the reflection $\refl$.
	The square is one point of an antipodeal
	pair invariant under reflection in the diameter.}
\end{center}

It is easy to verify that the configurations are homometric 
without doing any calculations as follows.
Since $C_1 \cap C_2 = X \sqcup Y$ 
it suffices to show that the pairs 
$X \sqcup \{ p_1\}, \, X \sqcup \{ p_2\}$
and
$Y \sqcup \{ p_1\}, \, Y \sqcup \{ p_2\}$
are isometric (and so homometric).
We have
\begin{eqnarray*}
\ant( X \sqcup \{ p_1\} ) & = &  X \sqcup \{ \ant(p_1)\}\\
\refl( Y \sqcup \{ p_1\} ) & = &  Y \sqcup \{ \refl(p_1)\}
\end{eqnarray*}

\subsection{Deforming the Mallows-Clark pair}

Using the decomposition it is easy to construct
a one parameter family of NCHP as follows.

The antipodeal map commutes with every isometry of the 
circle and in particular with the one parameter group of rotations $R_t$.
In particular, if $X$ is invariant under the antipodeal map
then so is $R_t(X)$.
For $t\in \RR$ define
\begin{eqnarray*}
C_1^t & :=& R_t(X) \sqcup Y \sqcup \{ p_1\}\\
C_2^t & :=& R_t(X) \sqcup Y \sqcup \{ p_2\}.
\end{eqnarray*}

Thus we have a family such that,
 $C_1^0$ and $C_2^0$ is the original Mallows-Clark pair
 and for all $t\in \RR$ 
 the configurations $C_1^t$ and $C_2^t$ 
 are homometric and non conjugate if $R_t(X)$
 is not invariant under $\refl$.

\subsection{Construction of homometric pairs via pairs of involutions }

We can now give a quite general construction for 
homometric pairs of subsets $C_1,C_2$ of a metric space $A$. 
It is interesting to note that  this construction
can be applied to Euclidean space 
to obtain pairs such that the Minkowski difference is not the same, that is
$$ C_1 - C_1 \neq  C_2 - C_2.$$
where the difference is defined by 
$$X - X := \{ x - y,\, x,y \in X \}.$$

\begin{thm}\label{struct_thm}
Let $A$ be a metric space  
and $\alpha_1,\alpha_2$ a pair of involutions
of $A$ acting by isometry.

If 
$X_1,X_2,P_1,P_2$ 
are  disjoint subsets of $A$
such that for $i = 1,2$
\begin{eqnarray}
\alpha_i(X_i) & = & X_i \\
\alpha_i(P_1) & = & P_2 
\end{eqnarray}

then the sets 
$$C_i := (X_1\sqcup X_2) \sqcup P_i,\, i = 1,2$$
are homometric. 

Further if  
\begin{enumerate}
\item $\alpha_1$ is the restriction of a unique isometry of $A$.
\item $X_1$ is the maximal subset of $C_1$ invariant under $\alpha_1$
\item $X_1$ is the unique  subset of  $C_1 \cup C_2$ conjugate  to  $X_1$
\end{enumerate}
then 
$C_1$ and $C_2$ are not conjugate.
\end{thm}

The conditions in the second part of the theorem  are  sufficient but not necessary.   
They hold for  original Mallows Clark pair where 
$\alpha_1$ is the antipodeal map and $X_1$ the antipodeal pair,
that is the unique set of points for which the
diameter of $C_i$ is attained.
So $X_1$ is invariant under any 
conjugation of $C_1$ and $C_2$.

In Paragraph \ref{triples} we will construct a pair which while 
not conjugate does not satisfy these conditions.
That is, there is $X_1$ 
conjugate to different  subsets of $X_1 \cup P_1$.
We will exploit this to produce a triple of homometric
pairwise non conjugate configurations.

\begin{figure}
\centering
\begin{subfigure}{.5\textwidth}
  \centering
 \includegraphics[scale=.5]{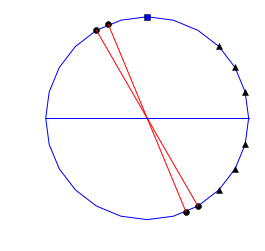}
\end{subfigure}%
\begin{subfigure}{.5\textwidth}
  \centering
  \includegraphics[scale=.5 ]{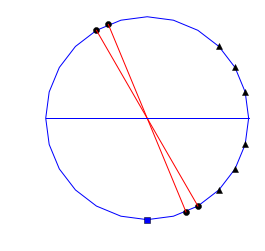}
\end{subfigure}
	\caption{The configurations of a generalized
	Mallows-Clark pair. 
	The black dots are antipodeal pairs. 
	The black triangles points of a set  invariant
	under reflection in the horizontal diameter.
	The squares are a pair of  antipodeal points
	invariant reflection in the diameter.}
\label{fig:test}
\end{figure}

\proof
The proof of homometry  is almost exactly as in Paragraph \ref{structure}.
Since $C_1 \cap C_2 = X_1 \sqcup X_2$ 
it suffices to show that, for $i = 1,2$, 
 the pairs 
$X_i \sqcup P_1 , \, X_i \sqcup  P_2$
are isometric (and so homometric).
Consider, as before,
\begin{eqnarray*}
\alpha_i( X_i \sqcup P_1 ) & = &  X_i \sqcup \alpha_i(P_1)
\end{eqnarray*}
for $i = 1,2$,
and by hypothesis $\alpha_i$ exchanges
 $P_1$ and $P_2$ so that 
$\alpha_i(P_1) = P_2$.

Now suppose that $C_1$ and $C_2$ are conjugate by an
(necessarily non trivial)  isometry  $\alpha$.
This isometry conjugates $X_1$ to a subset of $C_2$
so, by hypothesis, $\alpha(X_1) = X_1$ and,
since $\alpha_1$ is the restriction of   a unique isometry,
$$\alpha(X_2 \cup P_1) 
= \alpha_1(X_2 \cup P_1) = \alpha_1(X_2) \cup P_2$$
so that $\alpha_1(X_2) = X_2$
contradicting the maximality of $X_1$.

\hfill $\Box$

%
%
%
%
%

\subsection{Two constructions of (generalized) Mallows-Clark pairs}

We now use  Theorem \ref{struct_thm}
to construct the Mallows-Clark pair 
in two different ways.
Let $A$ be the circle $\{ |z| = 1, z \in \CC \}$
with the induced metric.
Now we make  two different choices for the involutions
\begin{enumerate}
\item using the antipodeal map and a reflection
$\alpha_1 : z \mapsto -z,\, \alpha_2 : z \mapsto \bar{z}$;
\item 
using a pair of distinct reflections 
$\alpha_1 : z \mapsto -\bar{z} ,\, \alpha_2 : z \mapsto \bar{z}$ .
\end{enumerate}
Note that all three of these involutions are restrictions of
  involutions of  the complex plane.
For the original Mallows-Clark pair this yields
 two different families of NCHP
whose intersection consists of  those configurations
such that $P_1 \cup P_2$ lie on the pair of lines
$ y = \pm x$.

\begin{coro} \label{mc families}
For  either of the choices of involution above
the   pair of configurations  arising from the
construction in Theorem \ref{struct_thm}
lies in a real analytic family of homometric pairs.
For almost all values of the parameters the
resulting pairs are non conjugate.

\end{coro}
\proof (Sketch) Existence follows  from the theorem.
The  almost all values  part is a consequence
of the fact that the family is real analytic.
\hfill $\Box$

Since the antipodeal map has no fixed points 
one sees that,
for  any family  obtained from the construction of   Corollary \ref{mc families},
 $X_1$ must have an even number of elements
whilst there is no restriction on the parity of
either $X_2$ or the $P_i$.
With this observation it is easy to see that:

\begin{coro}
For any $n\geq 4$ there is 
a family of homometric pairs on the circle
and for almost all values of the parameters the
resulting pairs are non conjugate.
\end{coro}

\begin{figure}
\centering
\begin{subfigure}{.5\textwidth}
  \centering
 \includegraphics[scale=.5]{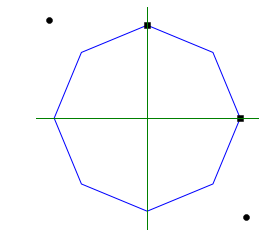}
\end{subfigure}%
\begin{subfigure}{.5\textwidth}
  \centering
  \includegraphics[scale=.5 ]{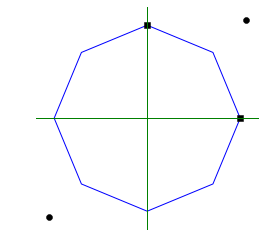}
\end{subfigure}
\caption{A (generalized)
	Mallows-Clark pair in the complex plane.
	}
\label{fig:pair2}
\end{figure}

\subsection{NCH Triples}
\label{triples}
\begin{figure}[h]
\centering
\begin{subfigure}{.5\textwidth}
  \centering
 \includegraphics[scale=.5]{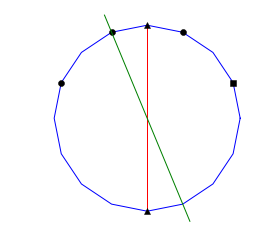}
\end{subfigure}%
\begin{subfigure}{.5\textwidth}
  \centering
  \includegraphics[scale=.5 ]{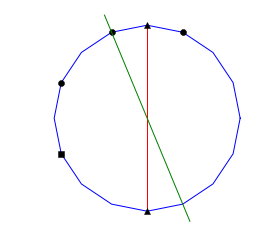}
\end{subfigure}
	\caption
	{A generalized
	Mallows-Clark pair associated a reflexion
	$\alpha_1$ and the antipodeal map $\alpha_2$.
	The invariant set of $\alpha_1$ consists of three black dots
	and $p_1$ is represented by a black square.
	}
\label{fig:triples stage1}
\end{figure}

\begin{figure}[h]
\centering
\begin{subfigure}{.5\textwidth}
  \centering
 \includegraphics[scale=.5]{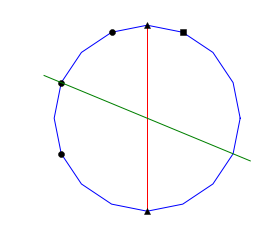}
\end{subfigure}%
\begin{subfigure}{.5\textwidth}
  \centering
  \includegraphics[scale=.5 ]{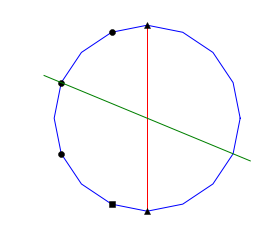}
\end{subfigure}
	\caption{
	The set $X_1$ is conjugate to a  subset 
	of $X_1'$ now represented by black dots in these figures.
	This new set is invariant under a different reflection
	allowing us to construct the third set of the triple.
	}
\label{fig:triples stage2}
\end{figure}

In this paragraph we construct a triple of homometric
pairwise non conjugate configurations.
Each pair of configurations
does not satisfy the conditions of Theorem \ref{struct_thm}:
the set  $X_1$  will be conjugate to different  subset $X_1'$
 of $X_1 \cup P_2$.
 The subset $X_1'$ is invariant under a reflexion $\alpha_1' \neq \alpha_1$
 and one defines $P_1'$ to be $ (X_1 \cup P_2) \setminus X_1'$.
Then using the construction of Theorem \ref{struct_thm}
one obtains 
a set  $$C_3 = X_1' \sqcup X_2  \sqcup \alpha_1'(P_1').$$
which  is homometric to $C_2 = X_1 \sqcup X_2 \sqcup P_2$.
This procedure is illustrated in Figure \ref{fig:triples stage1} and
Figure \ref{fig:triples stage2}.
Since the conditions of the second part of Theorem \ref{struct_thm}
no longer hold one must verify that the three configurations 
are pairwise non conjugate. 
First, observe that each of the pairwise intersections
of the  configurations is a unique antipodeal pair $X_2$
so that, if any pair is conjugate, 
the conjugation must be the antipodeal map $\alpha_2$.
Now it is easy to see that for the  antipodeal map 
$\alpha_2(C_i) \neq C_j$.

\section{Appendix: Golomb's Polynomial Method}
 
 For completeness we give an exposition of S. Golomb's approach (see e.g. \cite{Bekir-Golomb}) to NCHP via generating functions.
 
 \subsection{Factoring auto-correlation functions}
 Let $X$ and $Y$ be a NCHP. 
We suppose that $X$ consists of the points
$x_0 = 0 < x_1 < \ldots x_n = D,$
where $D$ is the diameter of $X$.
Golomb associates to $X$   the  polynomial of degree $D$
$$ r_X(t) :=  \sum_{x \in X} t^{d(x_0,x)}$$ 
It is easy to check that
$$r_X(t)r_X(1/t) =  \sum_{(x,y) \in X^2}  t^{x-y}.$$
Since $X$ is a subset of $\mathbb{R}$, $d(x,y) = |x-y|$,
and it follows immediately that the auto-correlation functions
$\sigma_X$ and $\sigma_Y$ are equal if and only if
\begin{equation} \label{equality}
r_X(t)r_X(1/t) = r_Y(t)r_Y(1/t).
\end{equation}
Golomb goes on to 
define 
$$r_X^*(t) :=  t^D r_{X}(1/t) = \sum_{x \in X} t^{d(x_n,x)}$$
and note  that $X$ is invariant under reflection if and only if
$r_X = r_X^*$.
Further, the condition  (\ref{equality})  is equivalent to
\begin{equation} \label{poly equality}
r_X r_X^* = r_Y r_Y^*.
\end{equation}

Golomb then made a remarkable observation:
the product $r_X(t)r^*_X(t)$ is a polynomial with integer coefficients  and  $\ZZ[t]$ is a unique factorisation domain (UFD)
so the  equation (\ref{poly equality})
tells us that any irreducible factor 
 of  $r_X$ is either a factor of $r_Y$ or $r_Y^*$.
 Consequently, one has a decomposition of the polynomials
 into products
 \begin{eqnarray*}
  r_X &=&
   \gcd(r_X,r_Y)  (\gcd(r_X,r_Y^*) / \gcd(r_Y,r_Y^*) ) \\
  r_ Y&=& \gcd(r_X,r_Y) 
  (  \gcd(r_Y,r_X^*)  /  \gcd(r_X,r_X^*) ).
  \end{eqnarray*}
 Note that $\gcd(r_X,r_Y^*) = \gcd(r_X^*,r_Y)$ so that 
 \begin{eqnarray} \label{decomposition}
 r_X = PQ,\, r_Y = PQ^*
 \end{eqnarray}
 where $P =    \gcd(r_X,r_Y)$ and $Q$ is the other factor.
 
 Since $\ZZ[t_1,t_2 \ldots t_n]$ is a UFD a similar 
 result is true in higher dimensions.

\subsection{An example} For example  for the NCHP $ \{0,1,4,10,12,17\}$ and $\{0,1,8,11,13,17\}$ cited in the introduction one has the 
 decomposition
\begin{eqnarray*}
r_X(t) = 1+t +t^4 +t^{10}+t^{12} +t^{17} 
&=&  (t^6 + t  + 1) (t^{11} - t^5 + t^4 + 1)\\
r_Y(t) = 1+t +t^8 +t^{11}+t^{13} +t^{17}
& = &  (t^6 + t  + 1)(1 - t^6 +t^7 +t^{11})
\end{eqnarray*}
 so that  $ \gcd(r_X,r_Y)  = t^6 + t  + 1$ 
 and the other factors form a pair $Q$ and $Q^*$.
 
 Using the decomposition  (\ref{decomposition})
 for $n \geq 9$ with at least two factors greater than or equal to $3$,
 one can construct NCHP with exactly $n$ points quite easily for example:
 \begin{eqnarray*}
(t^4 + t + 1) (t^7 + t^2 + 1)
&=& 
t^{11} + t^{8} + t^7 + t^{6} + t^{4} + t^{3} + t^{2} + t + 1\\
(t^4 + t^3 + 1) (t^7 + t^2 + 1) 
&=&
t^{11} + t^{10} + t^7 + t^6 + t^5 + t^4 + t^3 + t^2 + 1
 \end{eqnarray*}

 \subsection{Underlying geometric construction}
 
 There is an  underlying simple geometric construction
 which generalises to $\RR^n$ 
 (see \cite{AL} and Figures \ref{fig:test1} and \ref{fig:test2}  below).
 We say that a finite set $X \subset \RR$,
\textit{asymmetric} iff it is not invariant under
the inversion that swaps the minimum and maximum 
of $X$, namely
$$x \mapsto - x + \min(X)  + \max(X).$$
 Let $-X$ denote the set $\{ -x \mid x \in X \}$ as usual  
 and $\tau$  the  translation  $x\mapsto x + 1$.
Let 
$X_0 \subset ]-\frac12,\frac12[$ be a set
satisfying $ \min(X_0)  + \max(X_0) = 0$.
Note that the diameter of $X_0$ 
is in fact  $2\max(X_0)$.
Then for any  $X \subset \ZZ^+$
we can form the \textit{Minkowski sums}
\begin{eqnarray*}
X_0 \oplus X
&=&
 \{ a + b\mid (a,b) \in X_0 \times X \} \\
 (-X_0)  \oplus X
&=&
 \{ -a + b\mid (a,b) \in X_0  \times X \}.
\end{eqnarray*}
Under the hypothesis on $X_0$
\begin{eqnarray}\label{disjoint union}
X_0 \oplus X  = \bigsqcup_{n \in X} \tau^n(X_0),
\end{eqnarray}
we will say that each translate  $ \tau^n(X_0)$ is a \textit{cluster}
with center $n \in X$.
Observe now that for any pair of clusters 
$\tau^n( X_0),\,\tau^m( X_0)$ in $X_0 \oplus X$
the translation  $\rho: x \mapsto -x + m + n$
maps them to the  clusters to 
$\tau^m( -X_0),\,\tau^n(- X_0)$ in 
$(-X_0) \oplus X$.
Using this observation one can quite easily show  that 
the auto-correlation functions of   $X_0 \oplus X$ and  $(-X_0)  \oplus X$ 
are the same.

Moreover,
provided  both  $X_0$ and $X$ are asymmetric, 
these sets are not congruent.
To see this note that, under our hypothesis on $X_0$,  
$$m= \min( (-X_0) \oplus X) = \min( X_0 \oplus X)
 = \min(X_0) + \min(X) $$
Let $n_0 = \min(X)$ and consider
$\tau^{n_0}(X_0) \subset X_0 \oplus X$
and
$\tau^{n_0}(-X_0) \subset (-X_0) \oplus X$.
Then one has
$$\min(X_0) + \min(X)  \in \tau^{n_0}(X_0) \cap \tau^{n_0}(-X_0)$$
and, 
since the diameter of $X_0$ is smaller than 
the minimal distance between distinct  points of $X$,
it follows that
$$B_m\left(\frac12\right)  \cap (\pm X_0 \oplus X) = \tau^{n_0}(\pm X_0).$$
So that if $ (-X_0) \oplus X  =  X_0 \oplus X$
then $ \tau^{n_0}(X_0) =  \tau^{n_0}(- X_0)$
contradicting the hypothesis that $X_0$ was asymmetric.
A similar argument can be used to show that
 $ (-X_0) \oplus X $  and $X_0 \oplus X$
 cannot be related by an inversion.

 \begin{figure}
\centering
\begin{minipage}{.5\textwidth}
  \centering
  \includegraphics[width=1\linewidth]{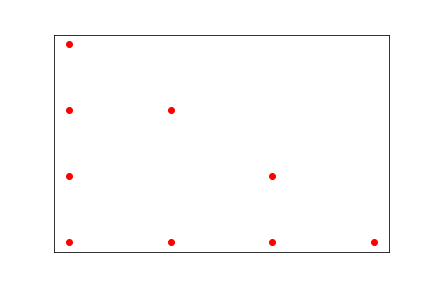}
  \captionof{figure}{$X_0 \oplus X$.}
  \label{fig:test1}
\end{minipage}%
\begin{minipage}{.5\textwidth}
  \centering
  \includegraphics[width=1\linewidth]{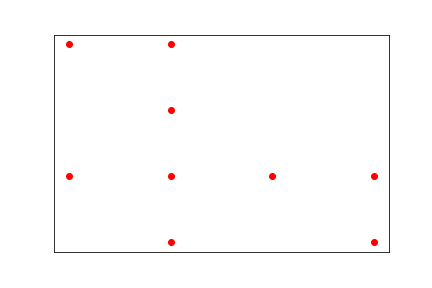}
  \captionof{figure}{$(-X_0) \oplus X$.}
  \label{fig:test2}
\end{minipage}
\end{figure}

Figures \ref{fig:test1} and \ref{fig:test2}  illustrate an example
of homometric sets  in $\RR^2$ obtained from
where $X_0 = \{ (0,0), (1,0) , (0,1)\}$ and $X = 2X_0$
by this process.

  \subsubsection{Nine point configurations} 
 Without the hypothesis on the diameter of $X_0$
 one  cannot guarantee that $X_0 \oplus X$ 
 decomposes as a disjoint union of translates as in (\ref{disjoint union}) above.
 So, in the general case, one has to deal with so-called 
 \textit{multi sets},
that is families of points counted with multiplicities,
  and the most convenient way to do 
 this appears to be via polynomials as above.
 
In order to satisfy the hypothesis that both $X_0$ and $X$
are asymmetric they must both have at least three points
and so the Minkowski sums each have at least nine points.
The construction of homometric configurations with less than nine points 
requires much more care as  can be seen 
from the factorisation in the previous section  
where one of the factors corresponds to a multiset where one point has a negative multiplicity
$$1+t +t^4 +t^{10}+t^{12} +t^{17} = (t^6 + t  + 1) (t^{11} - t^5 + t^4 + 1).$$

\end{document}